\documentclass[preprint,12pt]{elsarticle}

\usepackage{amssymb}
\usepackage{epstopdf}
\usepackage{color}

\newcommand{\udots}{\mathinner{\mskip1mu\raise1pt\vbox{\kern7pt\hbox{.}}
\mskip2mu\raise4pt\hbox{.}\mskip2mu\raise7pt\hbox{.}\mskip1mu}}
\usepackage{amsmath}

\newtheorem{theorem}{Theorem}
\newtheorem{corollary}{Corollary}
\newtheorem{lemma}{Lemma}

\begin{document}

\begin{frontmatter}



\title{ Doubly truncated moment risk measures for elliptical distributions
}

\author{Baishuai  Zuo}
\author{Chuancun Yin\corref{cor1}}
\cortext[cor1]{Corresponding author.}
\ead{ccyin@qfnu.edu.cn}

\address{School of Statistics and Data Science, Qufu Normal University, Qufu, Shandong 273165, P. R. China}

\begin{abstract}
 In this paper, we define doubly truncated moment (DTM), doubly truncated  skewness (DTS) and  kurtosis (DTK). We derive DTM formulae for elliptical family, with emphasis on normal, student-$t$, logistic, Laplace and Pearson type VII distributions. We also present explicit formulas of the DTE (doubly truncated expectation), DTV (doubly truncated variance), DTS and DTK for those distributions. As illustrative example, DTEs, DTVs, DTSs and DTKs of three industry segments' (Banks, Insurance, Financial and Credit Service)  stock return in London stock exchange are discussed.
\end{abstract}

\begin{keyword}

 Elliptical distributions;  Doubly truncated moment; Doubly truncated  skewness; Doubly truncated kurtosis; Normal distribution; Logistic distribution; Laplace distribution; Pearson type VII distribution
\end{keyword}

\end{frontmatter}

\baselineskip=20pt

\section{Introduction}
Landsman et al. (2016b) defined a new tail conditional moment (TCM) risk measure for a random variable $X$:
\begin{align}\label{(1)}
\mathrm{TCM}_{q}(X^{n})=\mathrm{E}\left[(X-\mathrm{TCE}_{q}(X))^{n}|X>x_{q}\right],
\end{align}
where
\begin{align*}
 \mathrm{TCE}_{q}(X)=E(X|X>x_{q})
 \end{align*}
 is tail conditional expectation (TCE) of $X$, $x_{q}$ is $q$-th quantile and $q\in(0,1)$.
Furthermore, they also defined novel types of tail conditional skewness and kurtosis (TCS and TCK):
\begin{align}\label{(2)}
\mathrm{TCS}_{q}(X)=\frac{\mathrm{E}\left[(X-\mathrm{TCE}_{q}(X))^{3}|X>x_{q}\right]}{\mathrm{TV}_{q}^{3/2}(X)}
\end{align}
and
\begin{align}\label{(3)}
\mathrm{TCK}_{q}(X)=\frac{\mathrm{E}\left[(X-\mathrm{TCE}_{q}(X))^{4}|X>x_{q}\right]}{\mathrm{TV}_{q}^{2}(X)}-3,
\end{align}
where
\begin{align*}
\mathrm{TV}_{q}(X)=\mathrm{E}[(X-\mathrm{TCE}_{X}(x_{q}))^{2}|X>x_{q}]
\end{align*}
is tail variance (TV) of $X$.
Since Landsman et al. (2016b) has been derived formulae of TCM for elliptical and log-elliptical distributions, and has been presented TCS and TCK for those distributions, Eini and Khaloozadeh (2021)  generalized those results to generalized skew-elliptical, Zuo and Yin (2021a) extended them to some shifted distributions.

Recently,  Roozegar et al. (2020) derived explicit expressions of the first two moments for doubly truncated multivariate normal mean-variance mixture distributions. Zuo and Yin (2021b) defined multivariate doubly truncated expectation and covariance risk measures, and derived formulas of multivariate doubly truncated expectation (MDTE) and covariance (MDTCov) of elliptical distributions. As special cases of  MDTE and MDTCov risk measures, authors also defined doubly truncated expectation (DTE) and variance (DTV) risk measures for a random variable $X$ as follows, respectively:
\begin{align}\label{(5)}
 \mathrm{DTE}_{(p,q)}(X)=\mathrm{E}(X|x_{p}<X<x_{q})
 \end{align}
 and
\begin{align}\label{(8)}
\mathrm{DTV}_{(p,q)}(X)=\mathrm{E}[(X-\mathrm{DTE}_{(p,q)}(X))^{2}|x_{p}<X<x_{q}],
\end{align}
 where $x_{k}$ $(k=p,~q)$ is $k$-th quantile, and $p,~q\in(0,1)$.

Inspired by those work, we define doubly truncated moment (DTM), and also define doubly truncated  skewness (DTS) and  kurtosis (DTK). Moreover, we derive doubly truncated moments (DTM) for elliptical family, and also give explicit expressions of DTE, DTV, DTS and DTK for this family and it's several special cases, such as  normal, student-$t$, logistic, Laplace and Pearson type VII distributions. As illustrative example, we discuss DTEs, DTVs, DTSs and DTKs of three industry segments' (Banks, Insurance, Financial and Credit Service)  stock return in London stock exchange.

The rest of the paper is organized as follows. Section 2 defines several doubly truncated risk measures. Section 3 introduces elliptical family and it's properties. In Section 4, we present $n$-th doubly truncated moments (TCM) for elliptical distributions, and derive explicit expressions of  DTV, DTS and DTK for this family. Special cases are given in Section 5. We give illustrative example in Section 6. Finally, in Section 7,  is the concluding remarks.
\section{Doubly truncated risk measures}
We define doubly truncated moment (DTM) risk measure of a random variable $X$ as follows:
\begin{align}\label{(4)}
\mathrm{DTM}_{(p,q)}(X^{n})=\mathrm{E}\left[(X-\mathrm{DTE}_{(p,q)}(X))^{n}|x_{p}<X<x_{q}\right],
\end{align}
where
$\mathrm{DTE}_{(p,q)}(X)$ is as in (\ref{(5)}), $x_{k}$ $(k=p,~q)$ is $k$-th quantile, and $p,~q\in(0,1)$.\\
 $\mathbf{Remark~1.}$ When $q\rightarrow1$, the doubly truncated moment (DTM) is reduced to tail conditional moment (TCM); When $p\rightarrow0$ and $q\rightarrow1$, the doubly truncated moment (DTM) is reduced to central moment.
Further, we define doubly truncated skewness (DTS) and kurtosis (DTK) risk measures:
\begin{align}\label{(6)}
\mathrm{DTS}_{(p,q)}(X)=\frac{\mathrm{E}\left[(X-\mathrm{DTE}_{(p,q)}(X))^{3}|x_{p}<X<x_{q}\right]}{\mathrm{DTV}_{(p,q)}^{3/2}(X)},
\end{align}
and
\begin{align}\label{(7)}
\mathrm{DTK}_{(p,q)}(X)=\frac{\mathrm{E}\left[(X-\mathrm{DTE}_{(p,q)}(X))^{4}|x_{p}<X<x_{q}\right]}{\mathrm{DTV}_{(p,q)}^{2}(X)}-3,
\end{align}
where
$\mathrm{DTV}_{(p,q)}(X)$ is as in (\ref{(8)}).\\
$\mathbf{Remark~2.}$ When $q\rightarrow1$, the doubly truncated skewness (DTS) is reduced to tail conditional skewness (TCS), and the doubly truncated kurtosis (DTK) is reduced to tail conditional kurtosis (TCK); When $p\rightarrow0$ and $q\rightarrow1$, the doubly truncated skewness (DTS) is reduced to skewness, and the doubly truncated kurtosis (DTK) is reduced to kurtosis.

Note that Molchanov and Cascos (2016), Cai et al. (2017) and Shushi and Yao (2020) proposed set risk measures (defined as a map from subset $S$ of
possible outcomes of losses $\Omega$ to some measure-valued space
$\mathcal{X}$, i.e., $S \subset \Omega \Rightarrow \rho(S) \in\mathcal{X}$) are
mathematically abstract and are very complicated when dealing
with risks. However,
tail conditional moment risk measures are relatively simpler
than that of the set risk measures and can be derived explicitly,
important for actuarial users (see Landsman et al., 2016a). In addition to these advantages of tail conditional moment, doubly truncated moment risk measures are more flexible than tail conditional moment risk measures. In other words, according to different needs, we can choose different $(p,q)$.
\section{Elliptical distributions}
A random variable $X$ is said to have an elliptically
symmetric distribution (see Landsman and Valdez, 2003)
\begin{align}\label{(9)}
f_{X}(x):=\frac{c_{1}}{\sigma}g_{1}\left\{\frac{1}{2}\left(\frac{x-\mu}{\sigma}\right)\right\},~x\in\mathbb{R},
\end{align}
 where $\mu$ is a location parameter, $\sigma>0$ is a scale parameter, $g_{1}(u)$, $u\geq0$, is the density generator of $X$, and denoted by $X\sim E_{1}(\mu,~\sigma^{2},~g_{1})$. The density generator $g_{1}$ satisfies the condition
\begin{align}\label{(10)}
\int_{0}^{\infty}s^{-1/2}g_{1}(s)\mathrm{d}s<\infty,
\end{align}
 and the normalizing constant  $c_1$ is given by
\begin{align*}
c_{1}&=\frac{\Gamma(1/2)}{(2\pi)^{1/2}}\left[\int_{0}^{\infty}s^{-1/2}g_{1}(s)\mathrm{d}s\right]^{-1}\\
&=\frac{1}{\sqrt{2}}\left[\int_{0}^{\infty}s^{-1/2}g_{1}(s)\mathrm{d}s\right]^{-1}.
\end{align*}
We define a sequence of cumulative generators $\overline{G}_{(k)},~k=1,~2\cdots,n,$:
\begin{align}\label{(11)}
\overline{G}_{(1)}(u)=\int_{u}^{\infty}g_{1}(s)\mathrm{d}s
\end{align}
and
\begin{align}\label{(11a)}
\overline{G}_{(k)}(u)=\int_{u}^{\infty}\overline{G}_{(k-1)}(s)\mathrm{d}s,~k\geq2.
\end{align}
The normalizing constants  $c_(k)^{\ast},~k\geq1,$ are given by
\begin{align}\label{(12)}
c_{(k)}^{\ast}=\frac{1}{\sqrt{2}}\left[\int_{0}^{\infty}s^{-1/2}\overline{G}_{(k)}(s)\mathrm{d}s\right]^{-1}.
\end{align}
 The density generators $\overline{G}_{(k)},~k=1,~2,\cdots,n,$ satisfy the condition
\begin{align}\label{(13)}
\int_{0}^{\infty}s^{-1/2}\overline{G}_{(k)}(s)\mathrm{d}s<\infty,~k=1,2\cdots,n.
\end{align}
\section{$N$-th doubly truncated moment}
In this section, we present $n$-th doubly truncated moment (DTM) of elliptical distributions, and also present DTV, DTS and DTK of elliptical distributions.

To derive $n$-th DTM of elliptical distributions, we define a new truncated distribution function as follows (see Zuo and Yin, 2021b): $$F_{Z}(a,b)=\int_{a}^{b}f_{Z}(z)\mathrm{d}z,$$
where $f_{Z}(z)$ is pdf of random variable $Z$.

 Firstly, we give following lemma.
\begin{lemma}\label{le.1}
Let $X\sim E_{1}(\mu,~\sigma^{2},~g_{1})$. Assume it satisfies conditions (\ref{(10)}) and (\ref{(13)}).
 Then
\begin{align}\label{(a9)}
\nonumber&\mathrm{E}[X^{n}|x_{p}<X<x_{q}]\\
&=\mu^{n}+n\mu^{n-1}\sigma\mathrm{DTE}_{(p,q)}(Y)+\sum_{i=2}^{n}\binom{n}{i}\mu^{n-i}\sigma^{i}\left[L_{1}+(i-1)\frac{c_{1}}{c_{(1)}^{\ast}}L_{2}\right],
\end{align}
where
$$\mathrm{DTE}_{(p,q)}(Y)=\frac{c_{1}\left(\overline{G}_{(1)}\left(\frac{1}{2}\xi_{p}^{2}\right)-\overline{G}_{(1)}\left(\frac{1}{2}\xi_{q}^{2}\right)\right)}{F_{Y}(\xi_{p},\xi_{q})},$$
$$L_{1}=\frac{c_{1}\left[\xi_{p}^{i-1}\overline{G}_{1}\left(\frac{1}{2}\xi_{p}^{2}\right)-\xi_{q}^{i-1}\overline{G}_{1}\left(\frac{1}{2}\xi_{q}^{2}\right)\right]}{F_{Y}(\xi_{p},\xi_{q})},$$
$$L_{2}=\frac{\int_{\xi_{p}}^{\xi_{q}}y^{i-2}c_{(1)}^{\ast}\overline{G}_{(1)}\left(\frac{1}{2}y^{2}\right)\mathrm{d}y}{F_{Y}(\xi_{p},\xi_{q})},$$
$\xi_{k}=\frac{x_{k}-\mu}{\sigma}$, $k=p,q$, and $Y\sim E_{1}(0,~1,~g_{1})$.
\end{lemma}
\noindent $\mathbf{Proof}$ Using definition, we have
\begin{align*}
 \mathrm{E}[X^{n}|x_{p}<X<x_{q}]=\frac{\int_{x_{p}}^{x_{q}}x^{n}\frac{c_{1}}{\sigma}g_{1}\left(\frac{1}{2}\left(\frac{x-\mu}{\sigma}\right)^{2}\right)\mathrm{d}x}{F_{X}(x_{p},x_{q})}.
\end{align*}
Applying the transformation $y=\frac{x-\mu}{\sigma}$, and using the Binomial Theorem, we obtain
\begin{align*}
 \mathrm{E}[X^{n}|x_{p}<X<x_{q}]&=\frac{\int_{\xi_{p}}^{\xi_{q}}(\sigma y+\mu)^{n}c_{1}g_{1}\left(\frac{1}{2}y^{2}\right)\mathrm{d}y}{F_{Y}(\xi_{p},\xi_{q})}\\
 &=\frac{\sum_{i=0}^{n}\binom{n}{i}\mu^{n-i}\sigma^{i}\int_{\xi_{p}}^{\xi_{q}} y^{i}c_{1}g_{1}\left(\frac{1}{2}y^{2}\right)\mathrm{d}y}{F_{Y}(\xi_{p},\xi_{q})}.
\end{align*}
Therefore,
\begin{align*}
&\mathrm{E}[X^{n}|x_{p}<X<x_{q}]\\
&=\mu^{n}+n\mu^{n-1}\sigma\mathrm{DTE}_{(p,q)}(Y)+\sum_{i=2}^{n}\binom{n}{i}\mu^{n-i}\sigma^{i}\left[L_{1}+(i-1)\frac{c_{1}}{c_{(1)}^{\ast}}L_{2}\right],
\end{align*}
as required.

Now we establish the formula of DTM for elliptical distributions.

\begin{theorem}\label{th.4} Suppose that $X\sim E_{1}(\mu,\sigma^{2},~g_{1})$, which satisfies conditions (\ref{(10)}) and (\ref{(13)}).
 Then
\begin{align}\label{(a10)}
 &\nonumber\mathrm{DTM}_{(p,q)}(X^{n})=(-1)^{n}\mathrm{DTE}_{(p,q)}^{n}(X)+(-1)^{n-1}n\mathrm{DTE}_{(p,q)}^{n}(X)\\
 &\nonumber~~~~~~~~~~~~~~~~~~+\sum_{k=2}^{n}\binom{n}{k}(-\mathrm{DTE}_{(p,q)}(X))^{n-k}[\mu^{k}+k\mu^{k-1}\sigma\mathrm{DTE}_{(p,q)}(Y)]\\
&+\sum_{k=2}^{n}\sum_{i=2}^{k}\binom{n}{k}\binom{k}{i}(-\mathrm{DTE}_{(p,q)}(X))^{n-k}\mu^{k-i}\sigma^{i}\left[L_{1}+(i-1)\frac{c_{1}}{c_{(1)}^{\ast}}L_{2}\right],~n\geq2,
\end{align}
where
\begin{align}\label{(a11)} \mathrm{DTE}_{(p,q)}(X)=\mu+\sigma\frac{c_{1}\left[\overline{G}_{1}\left(\frac{1}{2}\xi_{p}^{2}\right)-\overline{G}_{1}\left(\frac{1}{2}\xi_{q}^{2}\right)\right]}{F_{Y}(\xi_{p},\xi_{q})},
\end{align}
 $\xi_{k}$, $k=p,q$, $L_{1},$ $L_{2}$ and $Y$ are the same as those in Lemma 1.
\end{theorem}
\noindent $\mathbf{Proof}$ Using the Binomial Theorem and basic algebraic calculations, we have
\begin{align*}
\mathrm{DTM}_{(p,q)}(X^{n})&=\mathrm{E}[(X-\mathrm{DTE}_{(p,q)}(X))^{n}|x_{p}<X<x_{q}]\\
&=\mathrm{E}\left[\sum_{k=0}^{n}\binom{n}{k}X^{k}(-\mathrm{DTE}_{(p,q)}(X))^{n-k}|x_{p}<X<x_{q}\right]\\
&=\sum_{k=0}^{n}\binom{n}{k}(-\mathrm{DTE}_{(p,q)}(X))^{n-k}\mathrm{E}[X^{k}|x_{p}<X<x_{q}].
\end{align*}
By (45) of Zuo and Yin (2021b), we obtain that $\mathrm{DTE}_{(p,q)}(X)$ is as in (\ref{(a11)}).\\
Then, using Lemma 1, we obtain Eq.(\ref{(a10)}), as required.\\
$\mathbf{Remark~3.}$ When $q\rightarrow1$, the $n$-th tail conditional moment (TCM) for elliptical distribution is given by
\begin{align}\label{(a12)}
 \nonumber&\mathrm{TCM}_{p}(X^{n})=(-1)^{n}\mathrm{TCE}_{p}^{n}(X)+(-1)^{n-1}n\mathrm{TCE}_{p}^{n}(X)\\
&\nonumber~~~~~~~~~~~~~~~~+\sum_{k=2}^{n}\binom{n}{k}(-\mathrm{TCE}_{p}(X))^{n-k}[\mu^{k}+k\mu^{k-1}\sigma\mathrm{TCE}_{p}(Y)]\\
&+\sum_{k=2}^{n}\sum_{i=2}^{k}\binom{n}{k}\binom{k}{i}(-\mathrm{TCE}_{p}(X))^{n-k}\mu^{k-i}\sigma^{i}\left[L_{1}+(i-1)\frac{c_{1}}{c_{(1)}^{\ast}}L_{2}\right],~n\geq2,
\end{align}
where
\begin{align}\label{(a13)} \mathrm{TCE}_{p}(X)=\mu+\sigma\frac{c_{1}\overline{G}_{1}\left(\frac{1}{2}\xi_{p}^{2}\right)}{\overline{F}_{Y}(\xi_{p})},
\end{align}
$$L_{1}=\frac{c_{1}\xi_{p}^{i-1}\overline{G}_{1}\left(\frac{1}{2}\xi_{p}^{2}\right)}{\overline{F}_{Y}(\xi_{p})},$$
$$L_{2}=\frac{\int_{\xi_{p}}^{\infty}y^{i-2}c_{(1)}^{\ast}\overline{G}_{(1)}\left(\frac{1}{2}y^{2}\right)\mathrm{d}y}{\overline{F}_{Y}(\xi_{p})},$$
and $\overline{F}_{X}(\cdot)$ denotes tail function of $X$.\\
Note that (\ref{(a12)}) is the result of Theorem 1 in Landsman et al. (2016b).\\
$\mathbf{Remark~4.}$ Letting $p\rightarrow0$ and $q\rightarrow1$ in Theorem 1, the $n$-th central moment (CM) for elliptical distribution leads to
\begin{align}\label{(a14)}
 \nonumber\mathrm{CM}(X^{n})=&(-1)^{n}\mu^{n}+(-1)^{n-1}n\mu^{n}+\sum_{k=2}^{n}\binom{n}{k}(-\mu)^{n-k}\mu^{k}\\
&+\sum_{k=2}^{n}\sum_{i=2}^{k}\binom{n}{k}\binom{k}{i}(-1)^{n-k}\mu^{n-i}\sigma^{i}(i-1)\frac{c_{1}}{c_{(1)}^{\ast}}L_{2},~n\geq2,
\end{align}
where
$$L_{2}=\int_{-\infty}^{\infty}y^{i-2}c_{(1)}^{\ast}\overline{G}_{(1)}\left(\frac{1}{2}y^{2}\right)\mathrm{d}y.$$

Now, we give explicit expressions of DTV, DTS and DTK for elliptical distributions.
\begin{corollary}\label{co.1} Under conditions of Theorem 1, we have
\begin{align}\label{(a15)}
 &\mathrm{DTV}_{(p,q)}(X)=-\mathrm{DTE}_{(p,q)}^{2}(X)+\mu^{2}+2\mu\sigma\mathrm{DTE}_{(p,q)}(Y)+\sigma^{2}\left(L_{1}+\frac{c_{1}}{c_{(1)}^{\ast}}L_{2}\right),
\end{align}
where
\begin{align*}
 L_{1}=\frac{c_{1}\left[\xi_{p}\overline{G}_{(1)}\left(\frac{1}{2}\xi_{p}^{2}\right)-\xi_{q}\overline{G}_{(1)}\left(\frac{1}{2}\xi_{q}^{2}\right)\right]}{F_{Y}(\xi_{p},\xi_{q})},~
 L_{2}=\frac{F_{Y_{(1)}}(\xi_{p},\xi_{q})}{F_{Y}(\xi_{p},\xi_{q})},
\end{align*}
where
 $\mathrm{DTE}_{(p,q)}(X)$, $\xi_{k},~k=p,q,$ and $Y$ are the same as those in Theorem 1. In addition, $Y_{(1)}\sim E_{1}(0,~1,~\overline{G}_{(1)})$.
\end{corollary}
Note that (\ref{(a15)}) coincides with the result of (65) in Zuo and Yin (2021b). When $q\rightarrow1$, (\ref{(a15)}) is the result of (1.7) in  Furman and Landsman (2006).\\
\begin{corollary}\label{co.2} Under conditions of Theorem 1, we have
\begin{align}\label{(a16)}
 \nonumber&\mathrm{DTS}_{(p,q)}(X)\\
 \nonumber&=\mathrm{DTV}_{(p,q)}^{-3/2}(X)\bigg\{\sum_{k=2}^{3}\binom{3}{k}[-\mathrm{DTE}_{(p,q)}(X)]^{3-k}[\mu^{k}+k\mu^{k-1}\sigma\mathrm{DTE}_{(p,q)}(Y)]\\
 &+2\mathrm{DTE}_{(p,q)}^{3}(X)+3[\mu-\mathrm{DTE}_{(p,q)}(X)]\sigma^{2}\left(L_{1}+\frac{c_{1}}{c_{(1)}^{\ast}}L_{2}\right)+\sigma^{3}\left(L_{1}^{\ast}+2L_{2}^{\ast}\right)\bigg\},
\end{align}
where
\begin{align*}
 L_{1}^{\ast}=\frac{c_{1}\left[\xi_{p}^{2}\overline{G}_{(1)}\left(\frac{1}{2}\xi_{p}^{2}\right)-\xi_{q}^{2}\overline{G}_{(1)}\left(\frac{1}{2}\xi_{q}^{2}\right)\right]}{F_{Y}(\xi_{p},\xi_{q})},
\end{align*}
\begin{align*}
 L_{2}^{\ast}&=\frac{c_{1}\left[\overline{G}_{(2)}\left(\frac{1}{2}\xi_{p}^{2}\right)-\overline{G}_{(2)}\left(\frac{1}{2}\xi_{q}^{2}\right)\right]}{F_{Y}(\xi_{p},\xi_{q})},
\end{align*}
where
 $\mathrm{DTV}_{(p,q)}(X)$, $L_{1}$ and $L_{2}$ are the same as those in Corollary 1.
\end{corollary}
\begin{corollary}\label{co.3} Under conditions of Theorem 1, we have
\begin{align}\label{(a17)}
 \nonumber&\mathrm{DTK}_{(p,q)}(X)\\
 \nonumber&=\mathrm{DTV}_{(p,q)}^{-2}(X)\bigg\{-3\mathrm{DTE}_{(p,q)}^{4}(X)
 +6[\mu-\mathrm{DTE}_{(p,q)}(X)]^{2}\sigma^{2}\left(L_{1}+\frac{c_{1}}{c_{(1)}^{\ast}}L_{2}\right)\\
 &\nonumber~~+\sum_{k=2}^{4}\binom{4}{k}(-\mathrm{DTE}_{(p,q)}(X))^{4-k}[\mu^{k}+k\mu^{k-1}\sigma\mathrm{DTE}_{(p,q)}(Y)]\\
 &~~+4[\mu-\mathrm{DTE}_{(p,q)}(X)]\sigma^{3}\left(L_{1}^{\ast}+2L_{2}^{\ast}\right)+\sigma^{4}\left(L_{1}^{\ast\ast}+3L_{2}^{\ast\ast}\right)\bigg\}-3,
\end{align}
where
\begin{align*}
 L_{1}^{\ast\ast}=\frac{c_{1}\left[\xi_{p}^{3}\overline{G}_{(1)}\left(\frac{1}{2}\xi_{p}^{2}\right)-\xi_{q}^{3}\overline{G}_{(1)}\left(\frac{1}{2}\xi_{q}^{2}\right)\right]}{F_{Y}(\xi_{p},\xi_{q})},
\end{align*}
\begin{align*}
 L_{2}^{\ast\ast}&=\frac{c_{1}\left[\xi_{p}\overline{G}_{(2)}\left(\frac{1}{2}\xi_{p}^{2}\right)-\xi_{q}\overline{G}_{(2)}\left(\frac{1}{2}\xi_{q}^{2}\right)\right]}{F_{Y}(\xi_{p},\xi_{q})}+\frac{c_{1}}{c_{(2)}^{\ast}}\frac{F_{Y_{(2)}}(\xi_{p},\xi_{q})}{F_{Y}(\xi_{p},\xi_{q})}.
\end{align*}
Here
 $\mathrm{DTV}_{(p,q)}(X)$, $L_{1}^{\ast}$ and $L_{2}^{\ast}$ are the same as those in Corollary 2. In addition, $Y_{(2)}\sim E_{1}(0,~1,~\overline{G}_{(2)})$.
\end{corollary}

Note that (\ref{(a16)}) and (\ref{(a17)}) coincide with the results of (3.22) and (3.24) in Landsman et al. (2016b) as $q\rightarrow1$, repectively.

 \section{Special cases}

 In the following, we present DTV, DTS and DTK for several special members of univariate elliptical distributions, such as normal, student-$t$, logistic, Laplace and Pearson type VII distributions.\\
 $\mathbf{Example~1}$ (Normal distribution) Let $X\sim N_{1}(\mu,~\sigma^{2})$. In this case, the density generators are expressed:
 \begin{align*}
g_{1}(u)=\overline{G}_{(1)}(u)=\overline{G}_{(2)}(u)=\exp\{-u\},
\end{align*}
and the normalizing constants are written as:
\begin{align*}
c_{1}=c_{(1)}^{\ast}=c_{(2)}^{\ast}=(2\pi)^{-\frac{1}{2}}.
\end{align*}
 Then
\begin{align*}
 &\mathrm{DTV}_{(p,q)}(X)=-\mathrm{DTE}_{(p,q)}^{2}(X)+\mu^{2}+2\mu\sigma\mathrm{DTE}_{(p,q)}(Y)+\sigma^{2}\left(L_{1}+1\right),
\end{align*}
\begin{align*}
 \nonumber&\mathrm{DTS}_{(p,q)}(X)\\
 \nonumber&=\mathrm{DTV}_{(p,q)}^{-3/2}(X)\bigg\{\sum_{k=2}^{3}\binom{3}{k}[-\mathrm{DTE}_{(p,q)}(X)]^{3-k}[\mu^{k}+k\mu^{k-1}\sigma\mathrm{DTE}_{(p,q)}(Y)]\\
 &~~+2\mathrm{DTE}_{(p,q)}^{3}(X)+3[\mu-\mathrm{DTE}_{(p,q)}(X)]\sigma^{2}\left(L_{1}+1\right)+\sigma^{3}\left(L_{1}^{\ast}+2L_{2}^{\ast}\right)\bigg\},
\end{align*}
\begin{align*}
 \nonumber&\mathrm{DTK}_{(p,q)}(X)\\
 \nonumber&=\mathrm{DTV}_{(p,q)}^{-2}(X)\bigg\{-3\mathrm{DTE}_{(p,q)}^{4}(X)
 +6[\mu-\mathrm{DTE}_{(p,q)}(X)]^{2}\sigma^{2}\left(L_{1}+1\right)\\
 &\nonumber~~+\sum_{k=2}^{4}\binom{4}{k}(-\mathrm{DTE}_{(p,q)}(X))^{4-k}[\mu^{k}+k\mu^{k-1}\sigma\mathrm{DTE}_{(p,q)}(Y)]\\
 &~~+4[\mu-\mathrm{DTE}_{(p,q)}(X)]\sigma^{3}\left(L_{1}^{\ast}+2L_{2}^{\ast}\right)+\sigma^{4}\left[L_{1}^{\ast\ast}+3(L_{1}+1)\right]\bigg\}-3,
\end{align*}
where
\begin{align*}
\mathrm{DTE}_{(p,q)}(X)=\mu+\sigma\frac{\phi(\xi_{p})-\phi(\xi_{q})}{F_{Y}(\xi_{p},\xi_{q})},
\end{align*}
\begin{align*}
 L_{1}=\frac{\xi_{p}\phi(\xi_{p})-\xi_{q}\phi(\xi_{q})}{F_{Y}(\xi_{p},\xi_{q})},~
 L_{1}^{\ast}=\frac{\xi_{p}^{2}\phi(\xi_{p})-\xi_{q}^{2}\phi(\xi_{q})}{F_{Y}(\xi_{p},\xi_{q})},
\end{align*}
\begin{align*}
 L_{2}^{\ast}=\frac{\phi(\xi_{p})-\phi(\xi_{q})}{F_{Y}(\xi_{p},\xi_{q})},~
 L_{1}^{\ast\ast}=\frac{\xi_{p}^{3}\phi(\xi_{p})-\xi_{q}^{3}\phi(\xi_{q})}{F_{Y}(\xi_{p},\xi_{q})},
\end{align*}
$\xi_{k}=\frac{x_{k}-\mu}{\sigma}$, $k=p,q$, and $Y\sim N_{1}(0,~1)$. In addition, $\phi(\cdot)$ is pdf of $1$-dimensional standard normal distribution.\\
$\mathbf{Example~2}$ (Student-$t$ distribution). Let
 $X\sim St_{1}\left(\mu,~\sigma^{2},~m\right).$
 In this case, the density generators are expressed (for details see Zuo et al., 2021):
\begin{align*}
g_{1}(u)=\left(1+\frac{2u}{m}\right)^{-(m+1)/2},
\end{align*}
\begin{align*}
\overline{G}_{(1)}(u)=\frac{m}{m-1}\left(1+\frac{2u}{m}\right)^{-(m-1)/2}
\end{align*}
and
\begin{align*}
\overline{G}_{(2)}(u)=\frac{m^{2}}{(m-1)(m-3)}\left(1+\frac{2u}{m}\right)^{-(m-3)/2}.
\end{align*}
The normalizing constants are written as:
\begin{align*}
 c_{1}=\frac{\Gamma\left((m+1)/2\right)}{\Gamma(m/2)(m\pi)^{\frac{1}{2}}},
 \end{align*}
\begin{align*}
 \nonumber c_{(1)}^{\ast}&=\frac{(m-1)\Gamma(1/2)}{(2\pi)^{1/2}m}\left[\int_{0}^{\infty}u^{1/2-1}\left(1+\frac{2t}{m}\right)^{-(m-1)/2}\mathrm{d}u\right]^{-1}\\
 &=\frac{(m-1)}{m^{3/2}B(\frac{1}{2},~\frac{m-2}{2})},~if~m>2
 \end{align*} and
 \begin{align*}
 \nonumber c_{(2)}^{\ast}&=\frac{(m-1)(m-3)\Gamma(1/2)}{(2\pi)^{1/2}m^{2}}\left[\int_{0}^{\infty}u^{1/2-1}\left(1+\frac{2t}{m}\right)^{-(m-3)/2}\mathrm{d}u\right]^{-1}\\
 &=\frac{(m-1)(m-3)}{m^{5/2}B(\frac{1}{2},~\frac{m-4}{2})},~if~m>4,
 \end{align*}
 where $\Gamma(\cdot)$ and $B(\cdot,\cdot)$ are Gamma function and Beta function, respectively. Then
\begin{align*}
 \nonumber&\mathrm{DTV}_{(p,q)}(X)\\
 &=-\mathrm{DTE}_{(p,q)}^{2}(X)+\mu^{2}+2\mu\sigma\mathrm{DTE}_{(p,q)}(Y)+\sigma^{2}\left(L_{1}+\frac{m}{m-2}L_{2}\right),~m>2,
\end{align*}
\begin{align*}
 \nonumber&\mathrm{DTS}_{(p,q)}(X)\\
 \nonumber&=\mathrm{DTV}_{(p,q)}^{-3/2}(X)\bigg\{\sum_{k=2}^{3}\binom{3}{k}(-\mathrm{DTE}_{(p,q)}(X))^{3-k}[\mu^{k}+k\mu^{k-1}\sigma\mathrm{DTE}_{(p,q)}(Y)]\\
 \nonumber&+2\mathrm{DTE}_{(p,q)}^{3}(X)+3[\mu-\mathrm{DTE}_{(p,q)}(X)]\sigma^{2}\left(L_{1}+\frac{m}{m-2}L_{2}\right)+\sigma^{3}\left(L_{1}^{\ast}+2L_{2}^{\ast}\right)\bigg\},\\
 &~m>2,
\end{align*}
\begin{align*}
 \nonumber&\mathrm{DTK}_{(p,q)}(X)\\
 \nonumber&=\mathrm{DTV}_{(p,q)}^{-2}(X)\bigg\{-3\mathrm{DTE}_{(p,q)}^{4}(X)
 +6[\mu-\mathrm{DTE}_{(p,q)}(X)]^{2}\sigma^{2}\left(L_{1}+\frac{m}{m-2}L_{2}\right)\\
 &\nonumber~~+\sum_{k=2}^{4}\binom{4}{k}\mathrm{C}_{4}^{k}(-\mathrm{DTE}_{(p,q)}(X))^{4-k}[\mu^{k}+k\mu^{k-1}\sigma\mathrm{DTE}_{(p,q)}(Y)]\\
 &~~+4[\mu-\mathrm{DTE}_{(p,q)}(X)]\sigma^{3}\left(L_{1}^{\ast}+2L_{2}^{\ast}\right)+\sigma^{4}\left(L_{1}^{\ast\ast}+3L_{2}^{\ast\ast}\right)\bigg\}-3,~m>2,
\end{align*}
where
\begin{align*} \mathrm{DTE}_{(p,q)}(X)=\mu+\sigma\frac{\Gamma\left((m+1)/2\right)\sqrt{m}\left[\left(1+\frac{\xi_{p}^{2}}{m}\right)^{-(m-1)/2}-\left(1+\frac{\xi_{q}^{2}}{m}\right)^{-(m-1)/2}\right]}{\Gamma(m/2)(m-1)\sqrt{\pi}F_{Y}(\xi_{p},\xi_{q})},
\end{align*}
\begin{align*}
 L_{1}=\frac{\Gamma\left((m+1)/2\right)\sqrt{m}\left[\xi_{p}\left(1+\frac{\xi_{p}^{2}}{m}\right)^{-(m-1)/2}-\xi_{q}\left(1+\frac{\xi_{q}^{2}}{m}\right)^{-(m-1)/2}\right]}{\Gamma(m/2)(m-1)\sqrt{\pi}F_{Y}(\xi_{p},\xi_{q})},
\end{align*}
\begin{align*}
 L_{2}&=\frac{F_{Y_{(1)}}(\xi_{p},\xi_{q})}{F_{Y}(\xi_{p},\xi_{q})},
\end{align*}
\begin{align*}
 L_{1}^{\ast}=\frac{\Gamma\left((m+1)/2\right)\sqrt{m}\left[\xi_{p}^{2}\left(1+\frac{\xi_{p}^{2}}{m}\right)^{-(m-1)/2}-\xi_{q}^{2}\left(1+\frac{\xi_{q}^{2}}{m}\right)^{-(m-1)/2}\right]}{\Gamma(m/2)(m-1)\sqrt{\pi}F_{Y}(\xi_{p},\xi_{q})},
\end{align*}
\begin{align*}
 L_{2}^{\ast}&=\frac{\Gamma\left((m+1)/2\right)m^{3/2}\left[\left(1+\frac{\xi_{p}^{2}}{m}\right)^{-(m-3)/2}-\left(1+\frac{\xi_{q}^{2}}{m}\right)^{-(m-3)/2}\right]}{\Gamma(m/2)(m-1)(m-3)\sqrt{\pi}F_{Y}(\xi_{p},\xi_{q})},
\end{align*}
\begin{align*}
 L_{1}^{\ast\ast}=\frac{\Gamma\left((m+1)/2\right)\sqrt{m}\left[\xi_{p}^{3}\left(1+\frac{\xi_{p}^{2}}{m}\right)^{-(m-1)/2}-\xi_{q}^{3}\left(1+\frac{\xi_{q}^{2}}{m}\right)^{-(m-1)/2}\right]}{\Gamma(m/2)(m-1)\sqrt{\pi}F_{Y}(\xi_{p},\xi_{q})},
\end{align*}
\begin{align*}
 L_{2}^{\ast\ast}=&\frac{\Gamma\left((m+1)/2\right)m^{3/2}\left[\xi_{p}\left(1+\frac{\xi_{p}^{2}}{m}\right)^{-(m-3)/2}-\xi_{q}\left(1+\frac{\xi_{q}^{2}}{m}\right)^{-(m-3)/2}\right]}{\Gamma(m/2)(m-1)(m-3)\sqrt{\pi}F_{Y}(\xi_{p},\xi_{q})}\\
 &+\frac{m^{2}}{(m-2)(m-4)}\frac{F_{Y_{(2)}}(\xi_{p},\xi_{q})}{F_{Y}(\xi_{p},\xi_{q})},~m>4,
\end{align*}
$\xi_{k}=\frac{x_{k}-\mu}{\sigma}$, $k=p,q$, $Y\sim St_{1}(0,~1,~m)$, $Y_{(1)}\sim E_{1}(0,~1,~\overline{G}_{(1)})$ and $Y_{(2)}\sim E_{1}(0,~1,~\overline{G}_{(2)})$.\\
$\mathbf{Example~3}$ (Logistic distribution). Let $X\sim Lo_{1}\left(\mu,~\sigma^{2}\right)$. In this case, the density generators are expressed (for details see Zuo et al., 2021):
 \begin{align*}
g_{1}(u)=\frac{\exp(-u)}{[1+\exp(-u)]^{2}},
\end{align*}
\begin{align*}
\overline{G}_{(1)}(u)=\frac{\exp(-u)}{1+\exp(-u)}
\end{align*}
and
\begin{align*}
 \overline{G}_{(2)}(u)=\ln\left[1+\exp(-u)\right].
 \end{align*}
 The normalizing constants are written as:
 \begin{align*}
 c_{1}=\frac{1}{(2\pi)^{1/2}\Psi_{2}^{\ast}(-1,\frac{1}{2},1)},
\end{align*}
\begin{align*}
  c_{(1)}^{\ast}=\frac{1}{(2\pi)^{1/2}\Psi_{1}^{\ast}(-1,\frac{1}{2},1)}
\end{align*}
 and
 \begin{align*}
  c_{(2)}^{\ast}=\frac{1}{(2\pi)^{1/2}\Psi_{1}^{\ast}(-1,\frac{3}{2},1)}.
\end{align*} Then
\begin{align*}
 \nonumber&\mathrm{DTV}_{(p,q)}(X)\\
 &=-\mathrm{DTE}_{(p,q)}^{2}(X)+\mu^{2}+2\mu\sigma\mathrm{DTE}_{(p,q)}(Y)+\sigma^{2}\left[L_{1}+\frac{\Psi_{1}^{\ast}(-1,\frac{1}{2},1)}{\Psi_{2}^{\ast}(-1,\frac{1}{2},1)}L_{2}\right],
\end{align*}
\begin{align*}
 \nonumber&\mathrm{DTS}_{(p,q)}(X)\\
 \nonumber&=\mathrm{DTV}_{(p,q)}^{-3/2}(X)\bigg\{\sum_{k=2}^{3}\binom{3}{k}(-\mathrm{DTE}_{(p,q)}(X))^{3-k}[\mu^{k}+k\mu^{k-1}\sigma\mathrm{DTE}_{(p,q)}(Y)]\\
 &+2\mathrm{DTE}_{(p,q)}^{3}(X)+3[\mu-\mathrm{DTE}_{(p,q)}(X)]\sigma^{2}\left[L_{1}+\frac{\Psi_{1}^{\ast}(-1,\frac{1}{2},1)}{\Psi_{2}^{\ast}(-1,\frac{1}{2},1)}L_{2}\right]+\sigma^{3}\left(L_{1}^{\ast}+2L_{2}^{\ast}\right)\bigg\},
\end{align*}
\begin{align*}
 \nonumber&\mathrm{DTK}_{(p,q)}(X)=\\
 \nonumber&\mathrm{DTV}_{(p,q)}^{-2}(X)\bigg\{-3\mathrm{DTE}_{(p,q)}^{4}(X)
 +6[\mu-\mathrm{DTE}_{(p,q)}(X)]^{2}\sigma^{2}\left[L_{1}+\frac{\Psi_{1}^{\ast}(-1,\frac{1}{2},1)}{\Psi_{2}^{\ast}(-1,\frac{1}{2},1)}L_{2}\right]\\
 &\nonumber+\sum_{k=2}^{4}\binom{4}{k}(-\mathrm{DTE}_{(p,q)}(X))^{4-k}[\mu^{k}+k\mu^{k-1}\sigma\mathrm{DTE}_{(p,q)}(Y)]\\
 &+4[\mu-\mathrm{DTE}_{(p,q)}(X)]\sigma^{3}\left(L_{1}^{\ast}+2L_{2}^{\ast}\right)+\sigma^{4}\left(L_{1}^{\ast\ast}+3L_{2}^{\ast\ast}\right)\bigg\}-3,
\end{align*}
where
\begin{align*} \mathrm{DTE}_{(p,q)}(X)=\mu+\sigma\frac{\phi(\xi_{p})(1+\sqrt{2\pi}\phi(\xi_{q}))-\phi(\xi_{q})(1+\sqrt{2\pi}\phi(\xi_{p}))}{\Psi_{2}^{\ast}(-1,\frac{1}{2},1)F_{Y}(\xi_{p},\xi_{q})(1+\sqrt{2\pi}\phi(\xi_{p}))(1+\sqrt{2\pi}\phi(\xi_{q}))},
\end{align*}
\begin{align*}
 L_{1}=\frac{\xi_{p}\phi(\xi_{p})(1+\sqrt{2\pi}\phi(\xi_{q}))-\xi_{q}\phi(\xi_{q})(1+\sqrt{2\pi}\phi(\xi_{p}))}{\Psi_{2}^{\ast}(-1,\frac{1}{2},1)F_{Y}(\xi_{p},\xi_{q})(1+\sqrt{2\pi}\phi(\xi_{p}))(1+\sqrt{2\pi}\phi(\xi_{q}))},~
 L_{2}&=\frac{F_{Y_{(1)}}(\xi_{p},\xi_{q})}{F_{Y}(\xi_{p},\xi_{q})},
\end{align*}
\begin{align*}
 L_{1}^{\ast}=\frac{\xi_{p}^{2}\phi(\xi_{p})(1+\sqrt{2\pi}\phi(\xi_{q}))-\xi_{q}^{2}\phi(\xi_{q})(1+\sqrt{2\pi}\phi(\xi_{p}))}{\Psi_{2}^{\ast}(-1,\frac{1}{2},1)F_{Y}(\xi_{p},\xi_{q})(1+\sqrt{2\pi}\phi(\xi_{p}))(1+\sqrt{2\pi}\phi(\xi_{q}))},
\end{align*}
\begin{align*}
 L_{2}^{\ast}&=\frac{\ln(1+\sqrt{2\pi}\phi(\xi_{p}))-\ln(1+\sqrt{2\pi}\phi(\xi_{q}))}{\sqrt{2\pi}\Psi_{2}^{\ast}(-1,\frac{1}{2},1)F_{Y}(\xi_{p},\xi_{q})},
\end{align*}
\begin{align*}
 L_{1}^{\ast\ast}=\frac{\xi_{p}^{3}\phi(\xi_{p})(1+\sqrt{2\pi}\phi(\xi_{q}))-\xi_{q}^{3}\phi(\xi_{q})(1+\sqrt{2\pi}\phi(\xi_{p}))}{\Psi_{2}^{\ast}(-1,\frac{1}{2},1)F_{Y}(\xi_{p},\xi_{q})(1+\sqrt{2\pi}\phi(\xi_{p}))(1+\sqrt{2\pi}\phi(\xi_{q}))},
\end{align*}
\begin{align*}
 L_{2}^{\ast\ast}&=\frac{\xi_{p}\ln(1+\sqrt{2\pi}\phi(\xi_{p}))-\xi_{q}\ln(1+\sqrt{2\pi}\phi(\xi_{q}))}{\sqrt{2\pi}\Psi_{2}^{\ast}(-1,\frac{1}{2},1)F_{Y}(\xi_{p},\xi_{q})}+\frac{\Psi_{1}^{\ast}(-1,\frac{3}{2},1)}{\Psi_{2}^{\ast}(-1,\frac{1}{2},1)}\frac{F_{Y_{(2)}}(\xi_{p},\xi_{q})}{F_{Y}(\xi_{p},\xi_{q})},
\end{align*}
$\xi_{k}=\frac{x_{k}-\mu}{\sigma}$, $k=p,q$, $Y\sim Lo_{1}(0,~1)$, $Y_{(1)}\sim E_{1}(0,~1,~\overline{G}_{(1)})$ and $Y_{(2)}\sim E_{1}(0,~1,~\overline{G}_{(2)})$. \\
$\mathbf{Remark~5}$ Here $\Psi_{\kappa}^{\ast}(z,s,a)$ is the generalized Hurwitz-Lerch zeta function defined by (see
Lin et al., 2006)
$$\Psi_{\kappa}^{\ast}(z,s,a)=\frac{1}{\Gamma(\kappa)}\sum_{n=0}^{\infty}\frac{\Gamma(\kappa+n)}{n!}\frac{z^{n}}{(n+a)^{s}},$$
which has an integral representation
$$\Psi_{\kappa}^{\ast}(z,s,a)=\frac{1}{\Gamma(s)}\int_{0}^{\infty}\frac{t^{s-1}e^{-at}}{(1-ze^{-t})^{\kappa}}\mathrm{d}t,$$
where $\mathcal{R}(a)>0$, $\mathcal{R}(s)>0$ when $|z|\leq1~(z\neq1)$, $\mathcal{R}(s)>1$ when $z=1$.\\
 $\mathbf{Example~4}$ (Laplace distribution). Let $X\sim La_{1}\left(\mu,~\sigma^{2}\right)$. In this case, the density generators are expressed (for details see Zuo et al., 2021):
\begin{align*}
g_{1}(u)=\exp(-\sqrt{2u}),
\end{align*}
\begin{align*}
\overline{G}_{(1)}(u)=(1+\sqrt{2u})\exp(-\sqrt{2u})
\end{align*}
and
\begin{align*}
 \overline{G}_{(2)}(u)=(3+2u+3\sqrt{2u})\exp(-\sqrt{2u}).
 \end{align*}
The normalizing constants are written as:
\begin{align*}
c_{1}=\frac{1}{2},~c_{(1)}^{\ast}=\frac{1}{4},~ c_{(2)}^{\ast}=\frac{1}{16}.
 \end{align*}
 Then
\begin{align*}
 &\mathrm{DTV}_{(p,q)}(X)=-\mathrm{DTE}_{(p,q)}^{2}(X)+\mu^{2}+2\mu\sigma\mathrm{DTE}_{(p,q)}(Y)+\sigma^{2}\left(L_{1}+2L_{2}\right),
\end{align*}
\begin{align*}
 \nonumber&\mathrm{DTS}_{(p,q)}(X)\\
 \nonumber&=\mathrm{DTV}_{(p,q)}^{-3/2}(X)\bigg\{\sum_{k=2}^{3}\binom{3}{k}(-\mathrm{DTE}_{(p,q)}(X))^{3-k}[\mu^{k}+k\mu^{k-1}\sigma\mathrm{DTE}_{(p,q)}(Y)]\\
 &+2\mathrm{DTE}_{(p,q)}^{3}(X)+3[\mu-\mathrm{DTE}_{(p,q)}(X)]\sigma^{2}\left(L_{1}+2L_{2}\right)+\sigma^{3}\left(L_{1}^{\ast}+2L_{2}^{\ast}\right)\bigg\},
\end{align*}
\begin{align*}
 \nonumber&\mathrm{DTK}_{(p,q)}(X)\\
 \nonumber&=\mathrm{DTV}_{(p,q)}^{-2}(X)\bigg\{-3\mathrm{DTE}_{(p,q)}^{4}(X)
 +6[\mu-\mathrm{DTE}_{(p,q)}(X)]^{2}\sigma^{2}\left(L_{1}+2L_{2}\right)\\
 &\nonumber~~+\sum_{k=2}^{4}\binom{4}{k}(-\mathrm{DTE}_{(p,q)}(X))^{4-k}[\mu^{k}+k\mu^{k-1}\sigma\mathrm{DTE}_{(p,q)}(Y)]\\
 &~~+4[\mu-\mathrm{DTE}_{(p,q)}(X)]\sigma^{3}\left(L_{1}^{\ast}+2L_{2}^{\ast}\right)+\sigma^{4}\left(L_{1}^{\ast\ast}+3L_{2}^{\ast\ast}\right)\bigg\}-3,
\end{align*}
where
\begin{align*} \mathrm{DTE}_{(p,q)}(X)=\mu+\sigma\frac{(1+|\xi_{p}|)\exp(-|\xi_{p}|)-(1+|\xi_{q}|)\exp(-|\xi_{q}|)}{2F_{Y}(\xi_{p},\xi_{q})},
\end{align*}
\begin{align*}
 L_{1}=\frac{\xi_{p}(1+|\xi_{p}|)\exp(-|\xi_{p}|)-\xi_{q}(1+|\xi_{q}|)\exp(-|\xi_{q}|)}{2F_{Y}(\xi_{p},\xi_{q})},~
 L_{2}&=\frac{F_{Y_{(1)}}(\xi_{p},\xi_{q})}{F_{Y}(\xi_{p},\xi_{q})},
\end{align*}
\begin{align*}
 L_{1}^{\ast}=\frac{\xi_{p}^{2}(1+|\xi_{p}|)\exp(-|\xi_{p}|)-\xi_{q}^{2}(1+|\xi_{q}|)\exp(-|\xi_{q}|)}{2F_{Y}(\xi_{p},\xi_{q})},
\end{align*}
\begin{align*}
 L_{2}^{\ast}&=\frac{(3+\xi_{p}^{2}+3|\xi_{p}|)\exp(-|\xi_{p}|)-(3+\xi_{q}^{2}+3|\xi_{q}|)\exp(-|\xi_{q}|)}{2F_{Y}(\xi_{p},\xi_{q})},
\end{align*}
\begin{align*}
 L_{1}^{\ast\ast}=\frac{\xi_{p}^{3}(1+|\xi_{p}|)\exp(-|\xi_{p}|)-\xi_{q}^{3}(1+|\xi_{q}|)\exp(-|\xi_{q}|)}{2F_{Y}(\xi_{p},\xi_{q})},
\end{align*}
\begin{align*}
 L_{2}^{\ast\ast}&=\frac{\xi_{p}(3+\xi_{p}^{2}+3|\xi_{p}|)\exp(-|\xi_{p}|)-\xi_{q}(3+\xi_{q}^{2}+3|\xi_{q}|)\exp(-|\xi_{q}|)}{2F_{Y}(\xi_{p},\xi_{q})}+\frac{8F_{Y_{(2)}}(\xi_{p},\xi_{q})}{F_{Y}(\xi_{p},\xi_{q})},
\end{align*}
$\xi_{k}=\frac{x_{k}-\mu}{\sigma}$, $k=p,q$, $Y\sim La_{1}(0,~1)$, $Y_{(1)}\sim E_{1}(0,~1,~\overline{G}_{(1)})$ and $Y_{(2)}\sim E_{1}(0,~1,~\overline{G}_{(2)})$. In addition, $|\cdot|$ is absolute value function.\\
$\mathbf{Example~5}$ (Pearson type VII distribution). Let
$X\sim PVII_{1}\left(\mu,~\sigma^{2},~t\right).$
In this case, the density generators are expressed:
 \begin{align*}
 g_{1}(u)=(1+2u)^{-t},
 \end{align*}
 \begin{align*}
 \overline{G}_{(1)}(u)=\frac{1}{2(t-1)}(1+2u)^{-(t-1)}
 \end{align*}
 and
 \begin{align*}
 \overline{G}_{(2)}(u)=\frac{1}{4(t-1)(t-2)}(1+2u)^{-(t-2)}.
 \end{align*}
The normalizing constants are written as:
\begin{align*}
 c_{1}=\frac{\Gamma\left(t\right)}{\Gamma(t-1/2)\pi^{\frac{1}{2}}},~t>\frac{1}{2},
 \end{align*}
\begin{align*}
 c_{(1)}^{\ast}=\frac{2(t-1)}{B(\frac{1}{2},t-\frac{3}{2})},~t>\frac{3}{2}
 \end{align*}
 and
 \begin{align*}
 c_{(2)}^{\ast}=\frac{4(t-1)(t-2)}{B(\frac{1}{2},t-\frac{5}{2})},~t>\frac{5}{2}.
 \end{align*}
Then
\begin{align*}
 \nonumber&\mathrm{DTV}_{(p,q)}(X)\\
 &=-\mathrm{DTE}_{(p,q)}^{2}(X)+\mu^{2}+2\mu\sigma\mathrm{DTE}_{(p,q)}(Y)+\sigma^{2}\left(L_{1}+\frac{1}{2t-3}L_{2}\right),~t>\frac{3}{2},
\end{align*}
\begin{align*}
 \nonumber&\mathrm{DTS}_{(p,q)}(X)\\
 \nonumber&=\mathrm{DTV}_{(p,q)}^{-3/2}(X)\bigg\{\sum_{k=2}^{3}\binom{3}{k}(-\mathrm{DTE}_{(p,q)}(X))^{3-k}[\mu^{k}+k\mu^{k-1}\sigma\mathrm{DTE}_{(p,q)}(Y)]\\
 \nonumber&+2\mathrm{DTE}_{(p,q)}^{3}(X)+3[\mu-\mathrm{DTE}_{(p,q)}(X)]\sigma^{2}\left(L_{1}+\frac{1}{2t-3}L_{2}\right)+\sigma^{3}\left(L_{1}^{\ast}+2L_{2}^{\ast}\right)\bigg\},\\
 &~t>\frac{3}{2},
\end{align*}
\begin{align*}
 \nonumber&\mathrm{DTK}_{(p,q)}(X)\\
 \nonumber&=\mathrm{DTV}_{(p,q)}^{-2}(X)\bigg\{-3\mathrm{DTE}_{(p,q)}^{4}(X)
 +6[\mu-\mathrm{DTE}_{(p,q)}(X)]^{2}\sigma^{2}\left(L_{1}+\frac{1}{2t-3}L_{2}\right)\\
 &\nonumber~~+\sum_{k=2}^{4}\binom{4}{k}(-\mathrm{DTE}_{(p,q)}(X))^{4-k}[\mu^{k}+k\mu^{k-1}\sigma\mathrm{DTE}_{(p,q)}(Y)]\\
 &~~+4[\mu-\mathrm{DTE}_{(p,q)}(X)]\sigma^{3}\left(L_{1}^{\ast}+2L_{2}^{\ast}\right)+\sigma^{4}\left(L_{1}^{\ast\ast}+3L_{2}^{\ast\ast}\right)\bigg\}-3,~t>\frac{3}{2},
\end{align*}
where
\begin{align*} \mathrm{DTE}_{(p,q)}(X)=\mu+\sigma\frac{\Gamma(t-1)\left[\left(1+\xi_{p}^{2}\right)^{-(t-1)}-\left(1+\xi_{q}^{2}\right)^{-(t-1)}\right]}{2\Gamma(t-\frac{1}{2})\sqrt{\pi}F_{Y}(\xi_{p},\xi_{q})},
\end{align*}
\begin{align*}
 L_{1}=\frac{\Gamma(t-1)\left[\xi_{p}\left(1+\xi_{p}^{2}\right)^{-(t-1)}-\xi_{q}\left(1+\xi_{q}^{2}\right)^{-(t-1)}\right]}{2\Gamma(t-\frac{1}{2})\sqrt{\pi}F_{Y}(\xi_{p},\xi_{q})},~
 L_{2}&=\frac{F_{Y_{(1)}}(\xi_{p},\xi_{q})}{F_{Y}(\xi_{p},\xi_{q})},
\end{align*}
\begin{align*}
 L_{1}^{\ast}=\frac{\Gamma(t-1)\left[\xi_{p}^{2}\left(1+\xi_{p}^{2}\right)^{-(t-1)}-\xi_{q}^{2}\left(1+\xi_{q}^{2}\right)^{-(t-1)}\right]}{2\Gamma(t-\frac{1}{2})\sqrt{\pi}F_{Y}(\xi_{p},\xi_{q})},
\end{align*}
\begin{align*}
 L_{2}^{\ast}&=\frac{\Gamma(t-2)\left[\left(1+\xi_{p}^{2}\right)^{-(t-2)}-\left(1+\xi_{q}^{2}\right)^{-(t-2)}\right]}{4\Gamma(t-\frac{1}{2})\sqrt{\pi}F_{Y}(\xi_{p},\xi_{q})},
\end{align*}
\begin{align*}
 L_{1}^{\ast\ast}=\frac{\Gamma(t-1)\left[\xi_{p}^{3}\left(1+\xi_{p}^{2}\right)^{-(t-1)}-\xi_{q}^{3}\left(1+\xi_{q}^{2}\right)^{-(t-1)}\right]}{2\Gamma(t-\frac{1}{2})\sqrt{\pi}F_{Y}(\xi_{p},\xi_{q})},
\end{align*}
\begin{align*}
 L_{2}^{\ast\ast}=&\frac{\Gamma(t-2)\left[\xi_{p}\left(1+\xi_{p}^{2}\right)^{-(t-2)}-\xi_{q}\left(1+\xi_{q}^{2}\right)^{-(t-2)}\right]}{4\Gamma(t-\frac{1}{2})\sqrt{\pi}F_{Y}(\xi_{p},\xi_{q})}\\
 &+\frac{1}{(2t-5)(2t-3)}\frac{F_{Y_{(2)}}(\xi_{p},\xi_{q})}{F_{Y}(\xi_{p},\xi_{q})},~t>\frac{5}{2},
\end{align*}
$\xi_{k}=\frac{x_{k}-\mu}{\sigma}$, $k=p,q$, $Y\sim PVII_{1}(0,~1,~t)$, $Y_{(1)}\sim E_{1}(0,~1,~\overline{G}_{(1)})$ and $Y_{(2)}\sim E_{1}(0,~1,~\overline{G}_{(2)})$.
\section{Illustrative example}
We discuss DTE, DTV, DTS and DTK of three industry segments in finance
(Banks $X_{1}$, Insurance $X_{2}$, Financial and Credit Service $X_{3}$) collecting stock return data in London stock exchange from April 2013 to November 2019 (For the data, see
https://finance.yahoo.com/ and https://www.londonstockexchange.com/) in the finance sector
of the market by using the results of parameter estimates in Shushi and Yao (2020). Using multivariate normal distribution to fit data. We denote it by
$$\mathbf{X}=(X_{1},X_{2},X_{3})^{T}\sim N_{3}(\boldsymbol{\mu},\mathbf{\Sigma}).$$
Parameters are computed using maximum likelihood estimation:
\begin{align*}
 \boldsymbol{\mu}=10^{-3}\left(\begin{array}{ccccccccccc}
-1.140677\\
5.896240\\
2.107343
\end{array}
\right),
\mathbf{\Sigma}=10^{-4}\left(\begin{array}{ccccccccccc}
19.088935&12.503116&-3.720492\\
12.503116&20.268816&-3.162601\\
-3.720492&-3.162601&8.851913
\end{array}
\right).
\end{align*}

(i) Considering $p+q=1$, let $(p,q)=(0.05,0.95),~(0.10,0.90),~(0.15,0.85),$ $(0.20,0.80),~(0.25,0.75),~(0.30,0.70),$ results are presented in Table 1, Figures 1, 2 and 3.

Table 1 and Figures 1-3 show DTEs, DTVs, DTSs and DTKs of $X_{1}$ (Banks), $X_{2}$ (Insurance) and $X_{3}$ (Financial and Credit Service) for $(p,q)=(0.05,0.95),~(0.10,0.90),~(0.15,0.85),$ $(0.20,0.80),~(0.25,0.75),~(0.30,0.70),$ respectively. In Table 1 we see that DTEs of Banks $X_{1}$ for different $(p,q)$ are same, and equal to mean $\mu_{1}=-1.140677\times10^{-3}$; DTEs of Insurance $X_{2}$ for different $(p,q)$ are same, and equal to mean $\mu_{2}=5.896240\times10^{-3}$; DTEs of Financial and Credit Service $X_{3}$ for different $(p,q)$ are same, and equal to mean $\mu_{3}=2.107343\times10^{-3}$;
DTEs of Insurance $X_{2}$ are the greatest, and DTEs of Banks $X_{1}$ are the least.

 As we see in Figure 1, there is a clear difference among the DTVs of three industry segments in finance. In three industry segments of finance, no matter how $p$ changes, DTVs of Financial and Credit Service $X_{3}$ are the least, and  DTVs of Insurance $X_{2}$ are
the greatest.  The
dispersion of values of the DTVs between Insurance $X_{2}$ and Financial and Credit Service $X_{3}$ is the largest at $(p,q)=(0.05,0.95)$. However, DTVs of Banks $X_{1}$, Insurance $X_{2}$, Financial and Credit Service $X_{3}$ are decreasing with increase of $p(<0.5)$,  which means that the smaller the volatility of data with increase of $p(<0.5)$

From Figure 2, we see that however $p$ changes, DTSs of  Insurance $X_{2}$, Financial and Credit Service $X_{3}$ are $0$.  This implies that for normal distribution, value of DTS for symmetric interval is $0$. It indicates that distribution is no skewness on interval. This can also explain why the skewness of the normal distribution is $0$.

We observe from Figure 3 that DTKs of Banks $X_{1}$, Insurance $X_{2}$, Financial and Credit Service $X_{3}$ are decreasing with increase of $p$. Furthermore, In three industry segments of finance, no matter how $p$ changes, DTKs of Banks $X_{1}$, Insurance $X_{2}$, Financial and Credit Service $X_{3}$ are same, which means that the value of DTK is not affected by expectation $\mu_{k}$ and variance $\sigma_{k}$ $(k=1,~2, ~ 3)$.

(ii) Considering $p-q=0.65$, let $(p,q)=(0.05,0.70),~(0.10,0.75),~(0.15,0.80),$ $(0.20,0.85),~(0.25,0.90),~(0.30,0.95),$ results are shown in Figures 4, 5, 6 and 7.

Figures 4-7 show DTEs, DTVs, DTSs and DTKs of $X_{1}$ (Banks), $X_{2}$ (Insurance) and $X_{3}$ (Financial and Credit Service) for $(p,q)=(0.05,0.70),$~\\$(0.10,0.75),~(0.15,0.80),$ $(0.20,0.85),~(0.25,0.90),~(0.30,0.95),$ respectively.
From Figure 4, we observe that DTEs of Banks $X_{1}$, Insurance $X_{2}$, Financial and Credit Service $X_{3}$ are increasing with increase of $p$.  In three industry segments of finance, increase rate of DTEs of Financial and Credit Service $X_{3}$ are the least, and increase rates of DTEs of  Banks $X_{1}$ and Insurance $X_{2}$  are almost equality. DTE of Financial and Credit Service $X_{3}$ is
the greatest for $p= 0.05$, and the DTEs of Insurance $X_{2}$ are the greatest for $p> 0.05$. However, the DTEs of Banks $X_{1}$ are
the least for $p\leq 0.25$, and DTE of Financial and Credit Service $X_{3}$ is
the least for $p= 0.3$.

As we see in Figure 5, there is a clear difference among the DTVs of three industry segments in finance. In three industry segments of finance, no matter how $p$ changes, DTVs of Financial and Credit Service $X_{3}$ are the least, and DTVs of Insurance $X_{2}$ are
the greatest. However, DTVs of Banks $X_{1}$, Insurance $X_{2}$, Financial and Credit Service $X_{3}$ are firstly decreasing then increasing with increase of $p$. At middle point, DTVs of Banks $X_{1}$, Insurance $X_{2}$, Financial and Credit Service $X_{3}$ are least. That can be explained by the fact that the closer to the mean (center), the smaller the degree of center deviation.

We observe from 6 that  DTSs of Banks $X_{1}$, Insurance $X_{2}$, Financial and Credit Service $X_{3}$ are increasing with increase of $p$. In three industry segments of finance, no matter how $p$ changes, DTSs of Banks $X_{1}$, Insurance $X_{2}$, Financial and Credit Service $X_{3}$ are same, which indicates that the value of DTS is not affected by expectation $\mu_{k}$ and variance $\sigma_{k}$ ($k=1,~2$,~$3$). For $p\leq 0.15$, the value of DTS is negative, it indicates that distribution is left skewness on interval.
For $p\geq 0.20$, the value of DTS is positive, it indicates that distribution is right skewness on interval.

From Figure 7, we notice that DTKs of Banks $X_{1}$, Insurance $X_{2}$, Financial and Credit Service $X_{3}$ are firstly decreasing then increasing with increase of $p$. At middle point, DTKs of Banks $X_{1}$, Insurance $X_{2}$, Financial and Credit Service $X_{3}$ are least. In three industry segments of finance, no matter how $p$ changes, DTKs of Banks $X_{1}$, Insurance $X_{2}$, Financial and Credit Service $X_{3}$ are same, which also indicates that the value of DTK is not affected by expectation $\mu_{k}$ and variance $\sigma_{k}$ ($k=1,~2$,~$3$).

From the above Table and Figures, it can be concluded that the values of DTS and DTK are not affected by expectation and variance for normal distribution. In addition, choosing different $(p,q)$ will may have different results to policy decision.
 \section{Concluding remarks}
 In this paper, we give the DTM risk measure for elliptical family, which is a generalization of TCM in Landsman et al. (2016b). The TCE, TV, TCS and TCK risk measures are extended to DTE, DTV, DTS and DTK risk measures, respectively. There are many special cases for this elliptical family. We give several special cases, including normal, student-$t$, logistic, Laplace and Pearson type VII distributions. Finally, we discuss  DTEs, DTVs, DTSs and DTKs of three industry segments' (Banks, Insurance, Financial and Credit Service)  stock return in London stock exchange, and conclude that choosing different $(p,q)$ will may have different results to policy decision. Furthermore, in Eini and Khaloozadeh (2021), the authors derive formula of TCM for generalized skew-elliptical distributions. It will, therefore, be of interest to extend the results established here to the
generalized skew-elliptical distributions.
\section*{Acknowledgments}
\noindent  The research was supported by the National Natural Science Foundation of China (No. 12071251)
\section*{Conflicts of Interest}
\noindent The authors declare that they have no conflicts of interest.

\section*{References}
\bibliographystyle{model1-num-names}







\end{document}